\documentclass[12pt,oneside]{article}
\usepackage{amsmath,amssymb,amsfonts,amsthm}
\usepackage{color}

\newcommand{\C}{{\cal C}}

\newcommand{\prof}{\noindent \textit{\textbf{Proof.\:\:}}}

\def\o#1{\overline{#1}}

\def\Section#1{\vspace{30truept}\addtocounter{section}{1}\setcounter{thm}{0}\setcounter{equation}{0}
{\noindent\Large\bf\arabic{section}.~~#1}\par \vspace{12pt}}
\newtheorem{thm}{Theorem}[section]

\newtheorem{lem}[thm]{Lemma}
\newtheorem{prop}[thm]{Proposition}

\newtheorem{rem}[thm]{Remark}

\numberwithin{equation}{section}

\begin{document}

\title{\textbf{CARTAN CONNECTIONS ASSOCIATED TO A $\beta$-CONFORMAL CHANGE
IN FINSLER GEOMETRY}}
\author{\textbf{S. H. Abed}}
\date{}
\maketitle

% End of preamble and beginning of text.
% Produces the title.
\vspace{-1.15cm}

\begin{center}
{Department of Mathematics, Faculty of Science,\\
Cairo University,\\ Giza, Egypt\\

\smallskip
E-mail: sabed@frcu.eun.eg}
\end{center}

\vspace{1cm} \vspace{1cm} \smallskip

\noindent \textbf{Abstract.} On a Finsler manifold $(M,L)$, we
consider the change
$L\longrightarrow\overline{L}(x,y)=e^{\sigma(x)}L(x,y)+\beta (x,y)$,
which we call a $\beta$-conformal change. This change generalizes
various types of changes in Finsler geometry: conformal,
$C$-conformal, $h$-conformal, Randers and generalized Randers
changes. Under this change, we obtain an explicit expression
relating the Cartan connection associated to $(M,L)$ and the
transformed Cartan connection associated to $(M,\overline{L})$. We
also express some of the fundamental geometric objects (canonical
spray, nonlinear connection, torsion tensors, ...etc.) of
$(M,\overline{L})$ in terms of the corresponding objects of $(M,L)$.
We characterize the $\beta$-homothetic change and give necessary and
sufficient conditions for the vanishing of the difference tensor in
certain cases.\newline It is to be noted that many known results of
Shibata, Matsumoto, Hashiguchi and others are retrieved as special
cases from this work.
\bigskip



\bigskip
%%%%%%%%%%%%%%%%%%%%%% Section %%%%%%%%%%%%%%%%%%%%%%%%%%%%
\newpage

\Section{Introduction and Notations}

Let $(M,L)$ be a Finsler space, where $M$ is an $n$-dimensional
$C^{\infty}$  manifold and $L(x,y)$ is a Finsler metric function. If
$\sigma (x)$ is\ a
function in each coordinate neighborhood of $M$, the change $%
L(x,y)\longrightarrow e^{\sigma (x)}L(x,y)$ is called a conformal
change. This change was introduced by M. S. Kneblman \cite{r7} and
deeply investigated by many authors: \cite{r23}, \cite{r5},
\cite{r6}, ... etc. The change $L(x,y)\longrightarrow L(x,y)+\beta
(x,y)$, where $L(x,y)$ is a Riemannian metric function and
$\beta=b_{i}y^{i}$ is a $1$-form on $M$, is called a Randers change
after Randers who first introduced it in \cite{r80}. The geometric
properties of such a change have been studied in various works:
\cite{r58},
 \cite{r69}, \cite{r37}, \cite{r12} ,... etc.
Matsumoto \cite {r1} introduced the transformation
$L^{*}(x,y)=L(x,y)+\beta (x,y)$, where $L(x,y)$ is a Finsler metric
function, and named it a $\beta $-change. He obtained the
relationship between the Cartan connection coefficients of $(M,L)$
and those of $(M,L^{*})$. Since then, this change has been
investigated by many authors: Shibata \cite{r3}, Miron
\cite{r37},...etc.

A change generalizing all the above mentioned changes has been
introduced by Abed \cite{r60} in the form:\vspace{-0.3cm}
\begin{equation}
L(x,y)\longrightarrow \o L(x,y)=e^{\sigma (x)}L(x,y)+\beta (x,y),
\label{eq.1}\vspace{-0.2cm}
\end{equation}%
where $\sigma $ is a function of $x$ and $\beta (x,y)=b_{i}(x)y^{i}$
is a 1-form on $M$, and named a $\beta$-conformal change. \footnote{
In  \cite{r60}, we called the change (\ref{eq.1}) a
{\lq\lq{conformal $\beta $-change}\rq\rq}, but we think that the
name \linebreak {\lq\lq{$\beta$-conformal change}\rq\rq} is rather
the appropriate one for such a change. This is the name that we will
always employ.} In fact, when $\beta =0$, the change (\ref{eq.1})
reduces to a conformal change. When $\sigma =0$, it reduces to a
$\beta $-change if $L~$is a Finsler metric function and to a Rander
change if $L$ is a Riemannian metric function. In \cite{r60}, we
have established the relationships between some important tensors
associated with $(M,L)$ and the corresponding tensors associated
with $(M,\overline{L})$. We have also investigated some invariant
and $\sigma $- invariant tensors.

In this paper, we still consider the $\beta $-conformal change
(\ref{eq.1}). Under this change we obtain an explicit expression
relating the Cartan connection $C\Gamma $ associated to the Finsler
manifold $(M,L)$ and the transformed Cartan connection
$C\overline{\Gamma }$ associated to the Finsler manifold
$(M,\overline{L})$ (cf. Theorem A). This result generalizes various
results of  Hashiguchi \cite{r23}, Izumi (\cite{r5}, \cite{r6}),
 Matsumoto \cite{r1}, Shibata \cite{r3} and others \cite{r4},...,etc.
(cf. Remark 3.1). Having established this crucial relation, we draw
some consequences and conclusions from our fundamental theorem. We
relate the two canonical sprays $S^{r}$, $\overline{S}^{r}$ and also
the two Cartan nonlinear connections $N^{r}_{j}$,
$\overline{N}^{r}_{j}$. We get the relation between the torsion
tensors $C^{r}_{ij},\, P^{r}_{ij},\,R^{r}_{ij}$ and the
corresponding tensors
$\overline{C}^{r}_{ij},\,\overline{P}^{r}_{ij},\,\overline{R}^{r}_{ij}$.
We terminate the paper by two theorems (cf. Theorems B and C) which,
roughly speaking, characterize the $\beta$-homothetic change and the
vanishing of the difference tensor (between the two Cartan
connections) in different cases.

It should finally be noted that many known results are retrieved as
special cases from the obtained results as indicated in different
places of the work.\newline

 \noindent\textbf{NOTATIONS}. Throughout
the present paper, $(M,L)$ denotes an $n$-dimensional $C^{\infty}$
Finsler manifold,  $(x^{i})$ denote the coordinates of any arbitrary
point of the base manifold $M$ and $(y^{i})$ a supporting element at
the same point. We use the following notations:\newpage \noindent
$\partial _{i}$: partial differentiation with respect to
$x^{i}$,\newline $\dot{\partial}_{i}$: partial differentiation with
respect to $y^{i},$\newline
$g_{ij}:=\frac{1}{2}\dot{\partial}_{i}\dot{\partial}_{j}L^{2}=
\dot{\partial}_{i}\dot{\partial}_{j}E$: the fundamental metric
tensor,\newline
$L_{i}:=\dot{\partial}_{i}L=g_{ij}L^{j}:=g_{ij}\frac{y^{j}}{L}$: the
normalized supporting element,\newline
$h_{ij}:=L\dot{\partial}_{i}L_{j}=LL_{ij}=g_{ij}-L_{i}L_{j}$: the
angular metric tensor,\newline
$C_{ijk}:=\frac{1}{2}\dot{\partial}_{k}(g_{ij})=\frac{1}{2}\dot{\partial}_{i}\dot{\partial}_{j}\dot{\partial}_{k}E$,\newline
$C_{jk}^{i}:=g^{il}C_{ljk}$: the (h)hv-torsion tensor,\newline
$\Gamma _{jk}^{i}$: the coefficients of the Cartan connection
$C\Gamma$,\newline $N_{j}^{i}:=(y^{k}\Gamma _{jk}^{i})$: the
coefficients of the canonical nonlinear connection,\newline
$C\Gamma:=(\Gamma _{jk}^{i},N_{j}^{i},C_{jk}^{i})$: the Cartan
connection associated to $(M,L)$,\newline $\delta _{k}:=\partial
_{k}-N_{k}^{r}~\dot{\partial}_{r}$,\newline $X^{i}_{j|k}:=\delta
_{k}X_{j}^{i}+X_{j}^{r}\Gamma _{rk}^{i}-X_{r}^{i}\Gamma _{jk}^{r}$:
the horizontal or (h)-covariant derivative of $ X_{j}^{i}$.\newline
Contraction by $y^{i}$ will be denoted by the index $0$. For
example, we write $\Gamma _{0j}^{i}$ for $\Gamma _{kj}^{i}y^{k}$.

%\newpage %%%%%%%%%%%%%%%%%%%%%%%%%%%% Section %%%%%%%%%%%%%%%%%%%%%%%%%%%%%
%%%%%%%%%%%%%%%%%%%%%%%%%%%%%%%%%%%%%%%%%%%%%%%%%%%%%%%%%%%%%%%%%%%%%%%%%%%%

\vspace{30truept}\addtocounter{section}{1}\setcounter{thm}{0}%
\setcounter{equation}{0} {\noindent{\Large
\textbf{\arabic{section}.~~Basic tensors associated to a
$\beta$-conformal change}}}

\vspace{0.3cm} In this section we introduce some basic tensors
associated to a $\beta$-conformal change. Consider the
$\beta$-conformal change (\ref{eq.1}). The relation between the
associated normalized covariant supporting elements is given by:
\begin{equation}
_{~~}\overline{L}_{i}(x,y)=e^{\sigma (x)}L_{i}(x,y)+b_{i}(x).
\label{eq.2}
\end{equation}%
Consequently, if we write $\,L_{ij}:= \dot{\partial}_{j}L_{i}\,$,
  $\,L_{ijk}:= \dot{\partial}_{k}L_{ij}\,$, ... etc., we get

\begin{equation}
\overline{L}_{ij}(x,y)=e^{\sigma (x)}L_{ij}(x,y).  \label{eq.3}
\end{equation}%

\begin{equation}
\overline{L}_{ijk}(x,y)=e^{\sigma (x)}L_{ijk}(x,y). \label{eq.3a}
\end{equation}%

\begin{equation}
\overline{L}_{ijkl}(x,y)=e^{\sigma (x)}L_{ijkl}(x,y). \label{eq.3b}
\end{equation}%

It is to be noted that Equation (\ref{eq.3}) is equivalent to (\ref%
{eq.1}); it is then a characterization of $\beta $-conformal
changes.
It is clear that\ $L_{ij}=\frac{h_{ij}}{L}$ $~$is~$\sigma $%
-invariant \ (A tensor K is $\sigma $- invariant if $~\overline{K}%
(x,y)=e^{\sigma }K(x,y)$ under the $\beta$-conformal change
(\ref{eq.1})).
\begin{lem}
Under a $\beta$-conformal change $\overline{L}(x,y)=e^{\sigma
(x)}L(x,y)+\beta (x,y),~$the relation between the fundamental metric tensors $%
g_{ij} $ and $\overline{g}_{ij}$ is given by: $\ $%
\begin{equation*}
\overline{g}_{ij}=\,\,\,\tau (g_{ij}-L_{i}L_{j})+\,\overline{L}_{i}\,%
\overline{L}_{j}
\end{equation*}%
and the relation between the corresponding covariant components
$g^{ij}$\
and $\overline{g}^{ij}$ is given by:%
\begin{equation}
\overline{g}^{ij}=\,\,\,\tau ^{-1}g^{ij}+\mu \,l\,^{i\,}l\,^{j}\vspace{0.15cm%
}-\tau ^{-2}(_{\,\,}l\,^{i}b\,^{j}+l\,^{j}b^{i}),  \label{eq.4}
\end{equation}%
where $\,\,\mu =(e^{\sigma }Lb^{2}+\beta)/\overline{L}\tau
^{2}$,\,\,\,\,$\tau =e^{\sigma }\frac{\overline{L}}{L}$,\,\,\,
$b^{2}=b_{i}b^{i}$\, and \, $b^{i}=g^{ij}b_{j}$.
\end{lem}

\par
Differentiation the angular matric $h_{ij}$ with respect to $y^{k},$
we get

\begin{equation}
\dot{\partial}_{k}h_{ij}=2C_{ijk}-L^{-1}(L_{i}h_{jk}+L_{j}h_{ik}),
\label{eq.6}
\end{equation}%
from which we obtain

\begin{equation}
L_{ijk}=\frac{2}{L}C_{ijk}-\frac{1}{L^{2}}%
(h_{ij}L_{k}+h_{jk}L_{i}+h_{ki}L_{j}).  \label{eq.7}
\end{equation}%
Taking (\ref{eq.7}) into account, (\ref{eq.3a}) can be rewritten in terms of $%
C_{ijk}$ in the form

\begin{equation}
\overline{C}_{ijk}=\tau \lbrack
C_{ijk}+\frac{1}{2\overline{L}}h_{ijk}], \label{eq.8}
\end{equation}%
where $\,\,h_{ijk}=h_{ij}m_{k}+h_{jk}m_{i}+h_{ki}m_{j}$\, and\,
$\,\,m_{i}=b_{i}-\frac{\beta }{L}L_{i}$.\newline
\smallskip
(Note that $m_{0}=m_{i}y^{i}=0.$).

\bigskip
Using (\ref{eq.4}) and (\ref{eq.8}) we get the following

\begin{lem}
Under a $\beta $-conformal change $\overline{L}(x,y)=e^{\sigma
(x)}L(x,y)+\beta (x,y),$ the relation between the the (h)hv-torsion tensors$%
~C_{jk}^{i}~$and~$\overline{C}_{jk}^{i}~$has the form:
\begin{equation*}
\overline{C}_{jk}^{i}=C_{jk}^{i}+A_{jk}^{i},
\end{equation*}%
\vspace{-0.15cm} where
\begin{equation}
\
A_{jk}^{i}=\frac{1}{2\overline{L}}(h_{jk}\,m^{i}+h_{j}\,^{i}\,m_{k}+h^{i}
\,_{j}\,m_{k})-\,\frac{1}{\tau}
\,C_{jks}L^{i}\,b^{s}-\frac{1}{2\overline{L}\tau }
(2m_{j}\,m_{k}+m^{2}\,h_{jk})L^{i}, \label{eq.9}
\end{equation}
$m^{i}=g^{ik}\,m_{k}\,$ and $\,h_{j}\,^{i}:=g^{ik}\,h_{jk}$.
\end{lem}
\bigskip

Differentiation both sides of (\ref{eq.7}) with respect to $y^{h}$,
we get
\begin{align*}
L_{hijk}& =\frac{2}{L}\dot{\partial}_{k}C_{hij}-\frac{2}{L^{2}}%
(L_{h}C_{ijk}+L_{i}C_{hjk}+L_{j}C_{hik}+L_{k}C_{hij}) \\
& -\frac{1}{L^{3}}(h_{hi}h_{jk}+h_{hj}h_{ki}+h_{hk}h_{ij}) \\
& +\frac{2}{L^{3}}%
(h_{hi~}L_{j}L_{k}+h_{hj}L_{k}L_{i}+h_{hk}L_{i}L_{j}+h_{ij}L_{k}L_{h}+h_{jk}L_{i}L_{h}+h_{ki}L_{j}L_{h}).
\end{align*}%
Taking the above equation into account, (\ref{eq.3b}) gives the
relation
between $\dot{\partial}_{r}\overline{C}_{ijk}$ and $\dot{\partial}%
_{k}C_{hij} $:

\begin{lem}Under a $\beta $-conformal change $\overline{L}(x,y)=e^{\sigma
(x)}L(x,y)+\beta (x,y)$, we have
\begin{equation*}
\dot{\partial}_{r}\overline{C}_{ijk}=\tau \dot{\partial}_{r}C_{ijk}+\frac{%
e^{\sigma }}{L}C_{ijk}\,m_{r}+\mathfrak{S}_{i,j,k}[\frac{e^{\sigma }}{L}%
C_{ijr}\,m_{k}-\frac{e^{\sigma }}{2L^{2}}(h_{ij}(n_{rk}+\frac{\beta}{L}%
h_{rk})+h_{ir}\,n_{jk})],
\end{equation*}%
where $\ n_{rk}=m_{r}l_{k}+m_{k}l_{r}$ and $\mathfrak{S } _{i,j,k}$
denotes cyclic permutation on the indices $i,j,k$.
\end{lem}

%%%%%%%%%%%%%%%%%%%%%%%%%%%%%%  Section 3 %%%%%%%%%%%%%%%%%%%%%%%%%%%%%%%%%%%%%
%%%%%%%%%%%%%%%%%%%%%%%%%%%%%%%%%%%%%%%%%%%%%%%%%%%%%%%%%%%%%%%%%%%%%%%%%%%%%

\Section{Cartan connections associated to a
$\beta$-conformal\vspace{0.13cm} change}

This section is devoted to the determination of the relationship
between the Cartan
connection $C\Gamma $ associated to $(M,L)$ and the Cartan connection $C%
\overline{\Gamma }$ associated to $(M,\overline{L})$, under a $\beta $%
-conformal
change$\,\;L(x,y)\longrightarrow\overline{L}(x,y)=e^{\sigma(x)
}L(x,y)+\beta(x,y).$
\smallskip
\smallskip

Let $D_{jk}^{i}$ be the difference tensor between the Cartan
connection coefficients $\Gamma _{jk}^{i}~$and $\overline{\Gamma
}_{jk}^{i}$:
\begin{equation}
D_{jk}^{i}=\overline{\Gamma }_{jk}^{i}-\Gamma _{jk}^{i}
\label{eq.11}
\end{equation}%
Now, we are going to determine an explicit expression of
$D_{jk}^{i}.$ We will do so in three steps. Firstly, we determine
$D_{00}^{i}=D^{i}_{jk}\,y^{j}y^{k},~$then
$D_{0k}^{i}=D_{jk}^{i}\,y^{j}$ and finally $D_{jk}^{i}.$ Here is our
fundamental result:\newline \noindent

\noindent\textbf{Theorem A.} \emph{Under a $\beta$-conformal change
$\,\overline{L}(x,y)=e^{\sigma(x) }L(x,y)+\beta(x,y)$, the
relationship between the Cartan connection $\overline{\Gamma
}_{jk}^{r}\,$ of $\,(M,\overline{L})$ and the Cartan connection
$\Gamma_{jk}^{r}\,$ of $\,(M,L)$ is given by Equations (\ref{eq.40})
below.}

\vspace{7pt}
 \noindent \textit{\textbf{Proof.}} As we have said, we
proceed as follows: we firstly compute $D_{00}^{r}$, then
$D_{0j}^{r}$ and finally $D_{ij}^{r}.$ \vspace{0.25cm}
\par $\bullet$ Determination of
$D_{00}^{r}=D_{ij}^{r}\,y^{i}y^{j}$:
\newline
From (\ref{eq.3}), we get
\begin{equation}
\partial _{k}{}\overline{L}_{ij}=\partial _{k}(e^{\sigma }L_{ij}),
\label{eq.12}
\end{equation}
since
\begin{equation}
L_{ij|k}=\partial _{k}L_{ij}-L_{ijr}N_{k}^{r}-L_{rj}\Gamma
_{ik}^{r}-L_{ri}\Gamma _{jk}^{r}  \label{eq.13}
\end{equation}%
In virtue of $L_{ij\upharpoonleft k}=0$, (\ref{eq.13}) implies
\begin{equation}
\partial _{k}L_{ij}=L_{ijr}N_{k}^{r}+L_{rj}\Gamma _{ik}^{r}+L_{ri}\Gamma
_{jk}^{r}.  \label{eq.14}
\end{equation}%
Using (\ref{eq.14}) and (\ref{eq.11}), Equation (\ref{eq.12}) yields
\begin{equation}
L_{ijr}D_{0k}^{r}+L_{rj}D_{ik}^{r}+L_{ri}D_{jk}^{r}=\sigma
_{k}L_{ij}, \label{eq.15}
\end{equation}
where $\sigma _{k}=\frac{\partial \sigma }{\partial x^{k}}.$
\smallskip
\par
Now, differentiating (\ref{eq.2}) with respect to $x^{j},$ we get
\begin{equation}
\partial _{j}\overline{L}_{i}=\partial _{j}(e^{\sigma
}L_{i}+b_{i})=e^{\sigma }\sigma _{j}L_{i}+e^{\sigma }(\partial
_{j}L_{i})+\partial _{j}b_{i}.  \label{eq.16}
\end{equation}%
Taking into account the fact that $L_{i\upharpoonleft j}=0$ or,
equivalently, that $\partial _{j}L_{i}=L_{ri}N_{j}^{r}+L_{r}\Gamma
_{ji}^{r}~,$ we get
\begin{equation}
e^{\sigma }L_{ri}~D_{0j}^{r}+\overline{L}_{r}D_{ij}^{r}=e^{\sigma
}\sigma _{j}L_{i}+b_{i|j}.  \label{eq.17}
\end{equation}%
Equation (\ref{eq.17}) is equivalent to the following two equations:
\begin{equation}
e^{\sigma }(L_{ir}\text{\ }D_{0j}^{r}+L_{jr}\text{\ }D_{0i}^{r})+2\overline{L%
}_{r}D_{ij}^{r}=2E_{ij}+e^{\sigma }\sigma _{ij}.  \label{eq.19}
\end{equation}%
\begin{equation}
e^{\sigma }(L_{ir}\ D_{0j}^{r}-L_{jr}\ D_{0i}^{r})=2F_{ij}-e^{\sigma
}\mu _{ij},  \label{eq.20}
\end{equation}%
where \ \
\begin{eqnarray}
E_{ij} &=&\frac{1}{2}(b_{i|j}+b_{j|i}),\ \ \ \ \ F_{ij}=\frac{1}{2}%
(b_{i|j}-b_{j|i}),~  \notag \\
\ \ \sigma _{ij} &=&\sigma _{i}L_{j}+\sigma _{j}L_{i},\ \ \ \ \ \
\mu _{ij}=\sigma _{i}L_{j}-\sigma _{j}L_{i}.  \label{eq.21}
\end{eqnarray}%
$~$On the other hand, Equation (\ref{eq.15}) is equivalent to
\begin{equation}
2L_{jr}\
D_{ik}^{r}+L_{ijr}D_{0k}^{r}+L_{jkr}D_{0i}^{r}-L_{ikr}D_{0j}^{r}=\sigma
_{i}L_{jk}+\sigma _{k}L_{ij}-\sigma _{j}L_{ik}.  \label{eq.22}
\end{equation}%
Contracting (\ref{eq.19}) by $y^{j}$, we get

\begin{equation}
e^{\sigma }L_{ir}\ D_{00}^{r}+2\overline{L}_{r}D_{0i\text{ }%
}^{r}=2E_{0i}+e^{\sigma }(\sigma _{0}L_{i}+\sigma _{i}L),
\label{eq.23}
\end{equation}%
where we have put $\sigma _{0}:=\sigma _{j}y^{j}.$ Similarly, from (\ref%
{eq.20}) and (\ref{eq.22}), we get

\begin{equation}
e^{\sigma }L_{ir}\ D_{00}^{r}=2F_{i0}+e^{\sigma }(\sigma
_{0}L_{i}-\sigma _{i}L),  \label{eq.24}
\end{equation}

\begin{equation}
L_{ir}D_{0j}^{r}+L_{jr}\ D_{i0}^{r}+L_{ijr}D_{00}^{r}=\sigma
_{0}L_{ij}. \label{eq.25}
\end{equation}%
Again, contracting (\ref{eq.23}) by $y^{i}~~$gives

\begin{equation}
\overline{L}_{r}D_{00\text{ }}^{r}=E_{00}+e^{\sigma }\sigma _{0}L.
\label{eq.26}
\end{equation}%
Equations (\ref{eq.24}) and (\ref{eq.26}) can be written as a system
of algebraic equations in $D^{r}_{00}$:
\begin{eqnarray}
(i)~L_{ir}\ D_{00}^{r} &=&2e^{-\sigma }F_{i0}+(\sigma
_{0}L_{i}-\sigma
_{i}L)=B_{i},   \notag \\
(ii)~\overline{L}_{r}D_{00\text{ }}^{r} &=&E_{00}+e^{\sigma }\sigma
_{0}L=B. \label{eq.27}
\end{eqnarray}
\par
To determine the tensor $D^{r}_{00}$, we need the following lemma of
Matsumoto \cite{r1}:\newline \textit{\textbf{Lemma 1.}\ Given $B$\
and $B_{i}$ such that $B_{i}L^{i}=0$, the system of algebraic
equations\linebreak in $A^{r}$
\begin{equation*}
(i)\ L_{ri}\ A^{r}=B_{i}\ \ \ \ \ \ \ \ \ \ \ \ \ \ \ \ \ \ \ \
(ii)\ (L_{r}+b_{r})\ A^{r}\ =B,
\end{equation*}%
has the unique solution
\begin{equation*}
A^{r}~=LB^{r}~+\frac{L}{L+\beta }(B-LB_{\beta })\ L^{r},
\end{equation*}
where $B_{\beta }=B_{i}b^{i}.$}
\smallskip\par
Now, applying Lemma 1\, on the system (\ref{eq.27}), noting that $\,
B_{i}y^{i}=0$, we obtain an explicit expression for the required
tensor $D_{00}^{r}$:
\begin{equation}
D_{00}^{r}=2Le^{-\sigma
}F_{0}^{r}+\frac{L}{\overline{L}}(E_{00}-2Le^{- \sigma }F_{\beta
0})L^{r~}-L^{2}\sigma ^{r}+\frac{L}{\overline{L}} (2Le^{\sigma
}\sigma _{0}+L^{2}\sigma _{\beta })L^{r}, \label{eq.28}
\end{equation}
where \ $\sigma _{\beta }=\sigma _{i}\,b^{i}$,
$F_{0}^{r}=g^{ir}F_{i0}$\, and \,$F_{\beta 0}=F_{i0}\,b^{i}.$

\bigskip
\par $\bullet$ Determination of
$D_{0j}^{r}=D_{ij}^{r}\,y^{i}$:
\newline
Adding Equations (\ref{eq.20}) and (\ref{eq.25}), we get \noindent
\begin{equation}
L_{ir}\ D_{0j}^{r}=e^{-\sigma }F_{ij}-\frac{1}{2}L_{ijr}D_{00}^{r}+\frac{1}{2%
}(\sigma _{0}L_{ij}-\mu _{ij})=:G_{ij}.  \label{eq.29}
\end{equation}%
Equation (\ref{eq.23}) can be rewritten as

\begin{equation}
\overline{L}_{r}D_{i0}^{r}=E_{0i}-\frac{1}{2}e^{\sigma }L_{ir}D_{00}^{r}+%
\frac{1}{2}e^{\sigma }(\sigma _{0}L_{i}+\sigma _{i}L)=:G_{i}.
\label{eq.30}
\end{equation}
Substituting (\ref{eq.28}) in (\ref{eq.29}) and (\ref{eq.24}) in (\ref{eq.30}%
), we get the following expressions for $G_{ij}$ and $G_{i}$:
\begin{eqnarray}
(i)~G_{ij} &=&e^{-\sigma }F_{ij}-e^{-\sigma }LL_{ijr}F^{~r}{}_{0}+\frac{1}{2%
\overline{L}}L_{ij}(E_{00}-2Le^{-\sigma }F_{\beta 0})-  \notag \\
&&+\frac{1}{2}[L^{2}L_{ijr}\sigma ^{r}+\frac{L^{2}}{\overline{L}}%
L_{ij}~\sigma _{\beta }+(\sigma _{0}L_{ij}-\mu _{ij})]+\frac{L}{\overline{L}}%
e^{\sigma }\sigma _{0}L_{ij},  \notag \\
(ii)~G_{i} &=&E_{0i}-F_{0i}+e^{\sigma }\sigma _{i}~L~.\ \
\label{eq.31}
\end{eqnarray}%
From Equations (\ref{eq.27}), (\ref{eq.29}) and (\ref{eq.30}), the tensors $%
G_{ij}$ and $G_{j}$ have the properties:

\begin{equation}
G_{ij}y^{i}=0,\ \ \ \ \ \ \ \ G_{ij}y^{j}=B_{i},\ \ \ \ \ \ \ \
G_{j}y^{j}=B\ \label{eq.31a}
\end{equation}
\par To determine $D_{0j}^{r}$ we need the following lemma, which is
a generalized version of Matsumoto's lemma:\newline
\textit{\textbf{Lemma 2.}\, Given the tensor $G_{ij}$ and $G_{j}$,
with the properties (\ref{eq.31a}), the system of algebraic
equations in $D_{0j}^{r}$
\begin{equation}
(i)\ L_{ri}\ D_{0j}^{r}=G_{ij}\ \ \ \ \ \ \ \ \ \ \ \ \ \ \ \ \ \
(ii)\ \overline{L}_{r}\,D_{0j}^{r}=G_{j}\ .  \label{eq.32}
\end{equation}%
has the unique solution
\begin{equation}
D_{0j}^{i}=LG^{i}\,_{j}+\frac{L}{\overline{L}}(G_{j}-LG_{\beta
j})L^{i}, \label{eq.34}
\end{equation}%
where $~G^{i}\,_{j}=g^{ir}\,G_{rj}~$ and $~G_{\beta
j}=b^{i}\,G_{ij}.$ }\newline
\smallskip
\noindent Proof of Lemma 2.\, It follows, from the formula $\, g_{ij}=\,h_{ij}$ $%
+~L_{i}L_{j}$, that (i) of (\ref{eq.32}) can be rewritten as
\begin{equation}
{}g_{ir}D_{0j}^{r}=LG_{ij}+L_{i}L_{r}D_{0j}^{r}  \label{eq.35}
\end{equation}
Contaction of (\ref{eq.35}) by $b^{i}$ gives
\begin{equation}
b_{r}D_{0j}^{r}=LG_{\beta j}+\frac{\beta }{L}L_{r}D_{0j}^{r}
\label{eq.36}
\end{equation}%
On the other hand, taking (\ref{eq.36}) into account, (ii) of
(\ref{eq.32}) is rewritten as$\ \ \ \ \ \ \ \ \ \ \
\ \ \ \ \ \ \ \ \ \ $%
\begin{equation*}
e^{\sigma }L_{r}\ D_{0j}^{r}=G_{j}-b_{r}D_{0j}^{r}=G_{j}-[LG_{\beta
j}+\frac{\beta }{L}L_{r}D_{0j}^{r}]
\end{equation*}
which gives
\begin{equation}
L_{r}\ D_{0j}^{r}=\frac{L}{\overline{L}}\,[G_{j}-LG_{\beta j}]
\label{eq.37}
\end{equation}%
Substitution of (\ref{eq.37}) into (\ref{eq.35}) ends the proof of
the lemma.\newline The required tensor $D_{0j}^{r}$ is determined by
Equations (\ref{eq.34}).

\smallskip
\smallskip
\par $\bullet$ Determination of
$D_{ij}^{r}$:
\newline
Equations (\ref{eq.22}) and ({\ref{eq.19}) can be rewritten
respectively as}
\begin{eqnarray}
(i)~L_{ir}\ D_{jk}^{r}
&=&\frac{1}{2}\,(L_{jkr}D_{0i}^{r}-L_{ijr}D_{0k}^{r}-L_{ikr}D_{0j}^{r}+\sigma
_{j}L_{ik}+\sigma _{k}L_{ij}-\sigma _{i}L_{jk})=:H_{ijk}.\notag\\
(ii)~~\overline{L}_{r}D_{ij}^{r} &=&E_{ij}-\frac{1}{2}\,e^{\sigma }(L_{ir}\text{\ }%
D_{0j}^{r}+L_{jr}\text{\ }D_{0i}^{r}-\sigma _{ij})=:H_{ij}
\label{eq.38}
\end{eqnarray}%
From Equations (\ref{eq.32}) and (\ref{eq.38}), the tensors $H_{ijk}$ and $%
H_{jk}$ possess the properties:
\begin{equation}
H_{ijk}=H_{ikj},\ \ \ H_{ijk}y^{i}=0,\ \ \ H_{ijk}y^{j}=G_{ik},\
and\ \ \ H_{jk}y^{j}=G_{k}. \label{eq.38a}
\end{equation}
\smallskip
Finally, to determine $D_{jk}^{r}$ we need the next lemma, which
generalizes Lemma 2 and can be proved similarly:
\newline\smallskip\smallskip
\textit{\textbf{Lemma 3.}\, Given the tensor $H_{ijk}~$and $H_{jk}$,
with the properties (\ref{eq.38a}), the system of algebraic
equations in $D_{jk}^{r}$%
\begin{equation}
(i)\ L_{ir}\ D_{jk}^{r}=H_{ijk}\ \ \ \ \ \ \ \ \ \ \ \ \ \ \ \ \
(ii)\ \overline{L}_{r}\ D_{j~k}^{r}=H_{jk}.  \label{eq.39}
\end{equation}
has the unique solution
\begin{equation}
D_{~jk}^{i}=LH_{~jk}^{i}+\frac{L}{\overline{L}}(H_{jk}-LH_{\beta
jk})L^{i}, \label{eq.40}
\end{equation}
where $H_{~jk}^{i}=g^{ir}H_{rjk}$\ \ and \ $H_{\beta
jk}=b^{i}H_{ijk}.$}

\smallskip

Applying Lemma 3 on the system (\ref{eq.38}) having the properties
(\ref{eq.38a}), we obtain the expression (\ref{eq.40}) for the
required tensor $D_{~jk}^{i}$, which completes the proof of Theorem
3.1. \ $\Box$
\smallskip

It should be noted that the difference tensor \,$D_{ij}^{r}$\,
determined by (\ref{eq.40}) is explicitly expressed in terms of the
constituents of the $\beta$-conformal change (\ref{eq.1}) only.

\begin{rem}
Consider the $\beta$-conformal change (\ref{eq.1}):
$$L(x,y)\longrightarrow\overline{L}(x,y)=e^{\sigma(x)}L(x,y)+\beta(x,y).$$
\par $\bullet$ When the $\beta$-conformal change (\ref{eq.1}) is conformal
($\beta=0$), the difference tensor $D^{r}_{jk}$ takes the form:
\begin{eqnarray}
D^{r}_{jk}&=&L^{2}(C^{m}_{jk}\,C^{r}_{m}-C^{r}_{jm}\,C^{m}_{k}-C^{r}_{km}C^{m}_{j})\,
-(C^{r}_{jk}\,\sigma_{0}-C^{r}_{k}\,y_{j}-C^{r}_{j}\,y_{k}+C_{jk}\,y^{r}) \notag\\
&&+(\delta^{r}_{k}\,\sigma_{j}+\delta^{r}_{j}\,\sigma_{k}-g_{jk}\,\sigma^{r}),
\label{eq.40a}
\end{eqnarray}
where \ \ $C^{r}_{k}:=C^{r}_{kj}\,\sigma^{j}$ \,\, and \,\,
$C_{rk}:=C^{j}_{rk}\,\sigma_{j}$.\newline This is the case studied
by Hashiguchi \cite{r23}, Izumi (\cite{r5}, \cite{r6}) and others
\cite{r78}, \cite{r62}, ...etc.\newline
\par $\bullet$ When the $\beta$-conformal change (\ref{eq.1}) is $C$-conformal
($\beta=0$ and $C^{r}_{jk}\sigma^{j}=0$), the difference tensor
$D^{r}_{jk}$ takes the form:
\begin{equation}
D^{r}_{jk}=\delta^{r}_{k}\,\sigma_{j}+\delta^{r}_{j}\,\sigma_{k}-g_{jk}\,\sigma^{r}-C^{r}_{jk}\,\sigma_{0}.
\label{eq.41c} \end{equation} This is the case studied by Shibata
and Azuma \cite{r4}, Kim and Park
 \cite{r79} and others.\newline
\par $\bullet$ When the $\beta$-conformal change (\ref{eq.1}) is $h$-conformal
($\beta=0$ and
$C^{r}_{jk}\sigma_{r}=\frac{\C^{i}\sigma_{i}}{(n-1)}h_{jk}$), the
difference tensor $D^{r}_{jk}$ takes the form:\vspace{-0.3cm}
\begin{equation}
D^{r}_{jk}=\delta^{r}_{k}\,\rho_{j}+
\delta^{r}_{j}\,\rho_{k}-g_{jk}\,\rho^{r}-C^{r}_{jk}\,\rho-
\frac{C^{i}\sigma_{i}}{(n-1)}LL^{r}L_{j}L_{k},\vspace{-0.3cm}
\end{equation}
{\text{where }}\ \
$\rho_{j}:=\sigma_{j}+L\frac{\C^{i}\sigma_{i}}{(n-1)}L_{j}$\ \ \ and
$\rho=\rho_{j}y^{j}$.\newline This is the case studied by Izumi
 \cite{r6}.\newline
\par$\bullet$ When the
$\beta$-conformal change (\ref{eq.1}) is a Randers change
($\sigma=0$ and $L$ is Riemannian), the difference tensor
$D^{r}_{jk}$ takes the form:
\begin{eqnarray}
D_{~jk}^{i}&=&LH_{~jk}^{i}+\frac{L}{\overline{L}}(H_{jk}-LH_{\beta
jk})L^{i},\\
{\text{where}}\ \ \   H_{ij}&=&E_{ij}-\frac{1}{2}( G_{ij}+G_{ji}),\notag\\
H_{ijk}&=&\frac{1}{2}(L_{jkr}\,D^{r}_{0i}-L_{ijr}\,D^{r}_{0k}-L_{ikr}\,D^{r}_{0j}),\notag\\
D^{r}_{0i}&=&LG^{r}_{\ i}+\frac{L}{\overline{L}}\,(G_{i}-LG_{\beta
i})L^{r},\notag\label{eq.40b}
\end{eqnarray}
with $G_{ij}$ and $G_{i}$ given by Equation ({\ref{eq.41a}})
below.\newline This is the case studied by, Matsumoto \cite{r1},
Shibata, Shimada et al \cite{r12}  and others \cite{r9},...etc
\newline
\par$\bullet$ When the $\beta$-conformal change (\ref{eq.1}) is
a $\beta$-change ($\sigma=0$ and $L$ is Finslerian), the difference
tensor $D^{r}_{jk}$ takes the form (3.33), with $G_{ij}$ and $G_{i}$
given by Equation (\ref{eq.41b}) below.
\newline This is the case studied by Shibata \cite{r3}, Matsumoto
\cite{r1} and others.\newline
\par
The above discussion shows that our consideration is much more
general than various investigations existing in the literature.
\end{rem}

%%%%%%%%%%%%%%%%%%%%%%%%%%%%%%%%%%  Section 4 %%%%%%%%%%%%%%%%%%%%%%%%%%%%%%%%%%%%
%%%%%%%%%%%%%%%%%%%%%%%%%%%%%%%%%%%%%%%%%%%%%%%%%%%%%%%%%%%%%%%%%%%%%%%%%%%%%%%%%%

\Section{Consequences of the fundamental theorem}

Having obtained an explicit expression for the difference tensor
$D_{jk}^{r}$, we are in a position to draw some consequences and
conclusions from our fundamental theorem (Theorem A).

\begin{prop}
Under a $\beta$-conformal change, the following relations hold:
\begin{equation}
G_{ij}\ =\frac{1}{2}L_{ir}(\dot{\partial}_{j}~D_{00\text{ }}^{r})\ \
\ \ and\ \ \
G_{j}=\frac{1}{2}\overline{L}_{r}(\dot{\partial}_{j}~D_{00\text{
}}^{r})
\end{equation}
\end{prop}

\noindent \textit{\textbf{Proof.}} From equation (\ref{eq.24}), we
have
\begin{equation*}
\dot{\partial}_{j}(L_{ir}\ D_{00}^{r})=\dot{\partial}%
_{j}(2e^{-\sigma }F_{i0}+(\sigma _{0}L_{i}-\sigma _{i}L))
\end{equation*}%
which gives
\begin{equation*}
L_{ir}(\dot{\partial}_{j}~D_{00\text{ }}^{r})=2e^{-\sigma
}F_{ij~}-L_{ijr}D_{00}^{r}+(\sigma _{0}L_{ij}-\mu _{ij})=2G_{ij}.
\end{equation*}%
Similarly, the second relation is obtains from Equation
(\ref{eq.26}). \,$\Box $
\smallskip
\smallskip
\par
As a consequence of the above Proposition, taking (\ref{eq.32})\
into account, we get
\begin{equation*}
\dot{\partial}_{j}D_{00\text{ }}^{r}=2D_{0j\text{ }}^{r}
\end{equation*}

%\end{proof}

\begin{prop}
The tensor $G_{ij}$ can be written in a form free from the
difference
tensor $D^{r}_{ij}$\ as follows:%
\begin{eqnarray}
G_{ij} &=&e^{-\sigma }F_{ii~}+\frac{e^{-\sigma }}{L}
(L_{i}F_{j0}+L_{j}F_{i0})-2e^{-\sigma }C_{ij}^{k}\,F_{k0}+Gh_{ij}
\notag \\
&&+\frac{1}{2}[L^{2}L_{ijr}\sigma ^{r}+\frac{L^{2}}{\overline{L}}%
L_{ij}~\sigma _{\beta }+(\sigma _{0}L_{ij}-\mu
_{ij})]+\frac{L}{\overline{L}}\, e^{\sigma }\sigma _{0}L_{ij},
\label{eq.41}
\end{eqnarray}%
where we have put $\,G=\frac{1}{2L\overline{L}}(E_{00}-2Le^{-\sigma
}F_{\beta 0}).$
\end{prop}

In the case of a Randers metric, we retrieve the Matsumoto's tensor
\cite{r1}:
\begin{equation}
G_{ij}=F_{ij~}+\frac{1}{L}(L_{i}F_{j0}+L_{j}F_{i0})+Kh_{ij}\ \ \ and
\ \ \ G_{i}=E_{0i}-F_{0i},\label{eq.41a}
\end{equation}
where $\,\,\,\, K=\frac{1}{2L\overline{L}}(E_{00}-2LF_{\beta 0}).$
\bigskip
\par
In the case of a generalized Randers metric or a $\beta$-metric, we
get
\begin{equation}
G_{ij}=F_{ij~}+\frac{1}{L}(L_{i}F_{j0}+L_{j}F_{i0})-2C_{ij}^
{k}\,F_{k0}+Kh_{ij} \ \ \ and \ \ \
G_{i}=E_{0i}-F_{0i}.\label{eq.41b}
\end{equation}

\begin{prop}
Under~a~$\beta $-conformal change $L\longrightarrow
\,\overline{L}=e^{\sigma }L+\beta$, the canonical spray\,
$\overline{S}^{i}$\ and the Cartan nonlinear connection\,
$\overline{N}_{j}^{i}$ of $(M,\overline{L})$ can be determined in
terms of the corresponding objects of \,$(M,L)$\, in the form:
\begin{equation*}
\overline{S}^{r}\ =S^{r}+D_{00\text{ }}^{r},
\end{equation*}
\begin{equation}
\overline{N}_{j}^{r}~=N_{j}^{r}~+D_{0j\text{ }}^{r},  \label{eq.43}
\end{equation}
where $D_{00}^{r}$ and $D_{0j}^{i}$ are given by (\ref{eq.28} ) and
(\ref{eq.34}) respectively.
\end{prop}
\smallskip\par
It should be noted that the expression, relating the two nonlinear
connections, found by Hashiguchi \cite{r23} results as a special
case from the above proposition (by setting $\beta =0$).

\begin{prop}
Under a $\beta$-conformal change, the torsion tensors change as
follows:
\begin{description}
\item[(a)] $\overline{C}_{jk}^{i}=C_{jk}^{i}+A_{jk}^{i},$

\item[(b)] $\overline{P}_{jk}^{i}=P_{jk}^{i}-D_{jk}^{i}+B_{jk\text{ }}^{i},$

\item[(c)] $\overline{R}_{jk}^{i}=R_{jk}^{i}+\mathfrak{U}
_{j,k}(D_{0j|k}^{i}-(B_{jr\text{ }}^{i}+P_{jr}^{i})D_{0k}^{r}),$ \ \
\end{description}
where $A_{jk}^{i},~D_{0j}^{i}~$and $D_{jk}^{i}~$are given by
(\ref{eq.9}), (\ref{eq.34}) and (\ref{eq.40}) respectively,
$B_{jk}^{i}:=\dot\partial_{k}\,D_{0j}^{i}$  and $\,\mathfrak{U}
_{j,k}(Q_{jk}):=Q_{jk}-Q_{kj}$.
\end{prop}

\begin{prop}
The relation between the (v)hv-torsion tensors \ $\overline{P}_{hjk}$ and $%
P_{hjk}$ can be written in the form%
\begin{equation*}
\overline{P}_{hjk}=\tau P_{hjk}-\frac{\tau L}{2}%
[L_{hjkr}D_{00}^{r}+L_{hjr}D_{0k}^{r}+L_{hkr}D_{0j}^{r}+L_{jkr}D_{0h}^{r}-%
\sigma _{0}L_{hjk}].
\end{equation*}
\end{prop}

\noindent \textit{\textbf{Proof.}} By Proposition 4.4(b), taking
Equation (\ref{eq.40}) and the expressions
$B_{jk}^{i}:=\dot\partial_{k}\,D_{0j}^{i}$ and
$\overline{L}_{i}\ \overline{P}_{jk}^{i}=\overline{L}_{i}\ \overline{C}%
_{jk|0}^{i}=0$ into account, we get
\begin{equation*}
\ \ \ \ \overline{P}_{hjk}=\overline{g}_{ih}\
\overline{P}_{jk}^{i}=\,\,[\,\tau
(g_{ih}-L_{i}L_{h})+\,\overline{L}_{i}\,\overline{L}_{h}]~\overline{P}%
_{jk}^{i}
\end{equation*}
\begin{equation*}
\ \ \ \ \ \ =\tau
(g_{ih}-L_{i}L_{h})[P_{jk}^{i}+\dot{\partial}_{k}~D_{0j\text{
}}^{i}-D_{jk}^{i}]\
\end{equation*}
\begin{equation*}
=\tau P_{hjk}+\tau LL_{ih}(\dot{\partial}_{k}~D_{0j\text{ }}^{i}-D_{jk}^{i}).%
\newline
\end{equation*}%
\noindent Differentiation $L_{ir}\ D_{0j}^{r}~$with respect to
$y^{k}\ $and$\ $using
equations (\ref{eq.29}), we obtain%
\begin{equation*}
\dot{\partial}_{k}(L_{ir}\ D_{0j}^{r})=L_{irk}D_{0j}^{r}+L_{ir}\dot{\partial}%
_{k}D_{0j}^{r}\ \ \ \ \text{and \ \ \ }L_{ih}\dot{\partial}_{k}D_{0j}^{i}=%
\dot{\partial}_{k}G_{hj}-L_{ihk}D_{0j}^{h}\
\end{equation*}%
Therefore, $\overline{P}_{hjk}~$can be rewritten as

\begin{equation}
\overline{P}_{hjk}\ \ =\tau P_{hjk}+\tau L(\dot{\partial}%
_{k}G_{hj}-L_{hkr}D_{0j}^{r}\ -H_{hjk}).\label{eq.44}
\end{equation}
Now, by Proposition 4.1 and Equation (\ref{eq.29}), one can calculate $%
\dot{\partial}_{k}G_{hj}$ which will be of the form
\begin{equation*}
\dot{\partial}_{k}G_{hj}=-\frac{1}{2}\,L_{hjkr}D_{00}^{r}-L_{hjr}D_{0k}^{r}+\frac{1}{2}\,[
\sigma _{0}L_{hjk}+\sigma _{k}L_{hj}-\sigma _{h}L_{jk}+\sigma
_{j}L_{hk}].
\end{equation*}
Substituting the above equation into (\ref{eq.44}), we get the
result.$\,\,\Box $ \smallskip\smallskip
\par It shout be noted that
Proposition 4 of Matsumoto \cite{r1} results as a special case from
the above proposition (by setting $\sigma$=0 and letting $L$ be
Riemannian).\newline

\noindent\textbf{Theorem B.} \emph{Under~a~$\beta $-conformal change
$L\longrightarrow \,\overline{L}=e^{\sigma }L+\beta ,$ consider the
following two assertions:\newline {\em(i)} The covariant vector
$b_{i}~$is parallel with respect to the Cartan connection $C\Gamma
.$\newline {\em(ii)} The difference tensor $D_{jk}^{i}$vanishes
identically.\newline Then, we have\newline {\em(a)} If {\em(i)} and
{\em(ii)} hold, then $\sigma $~is homothetic.\newline {\em(b)} If
$\sigma $~is homothetic, then {\em(i)} and {\em(ii)} are
equivalent.}

\vspace{7pt}
 \prof
\par
(a) If $b_{i|j}=0,~$than $E_{ij}=0=F_{ij},~$by (\ref{eq.21}). This,
together with (\ref{eq.27}), imply
\begin{eqnarray*}
\sigma _{0}L_{i}-\sigma _{i}L &=&L_{ir}\ D_{00}^{r}, \\
~\ \ e^{\sigma }\sigma _{0}L &=&~\overline{L}_{r}D_{00\text{ }}^{r}.
\end{eqnarray*}%
Now, if, moreover, $D_{ij}^{r}=0,~$then $D_{00}^{r}=0,$ and so the
above two equations give
\begin{eqnarray*}
\sigma _{0}L_{i}-\sigma _{i}L &=&0,~ \\
\ \ \sigma _{0}L &=&0.
\end{eqnarray*}
These two equations, together, yield $\sigma _{i}=0~$and $\sigma
$~is homothetic.
\par
(b) Let $\sigma _{i}=0~$and $b_{i|j}=0.~$Then, it follows from (\ref
{eq.28}) that $D_{00}^{r}=0.$ This, together with (\ref{eq.29}) and
(\ref {eq.30}), imply that $G_{ij}=0~$and~$G_{i}=0.~$Consequently,
$D_{0j}^{r}=0$,  by (\ref{eq.34}). Again, this, together with\
(\ref{eq.38}), imply that $ H_{ijk}=0~ $and $H_{ij}=0.$
Consequently, $D_{ij}^{r}=0$, by (\ref{eq.40}).

On the other hand, let $\sigma _{i}=0$ and $D_{ij}^{r}=0.$ Then, (\ref{eq.29}%
) implies that $F_{ij}=0$ and (\ref{eq.38})(ii) implies that
$E_{ij}=0.~$ Consequently, it follows from (\ref{eq.21}) that
$b_{i|j}=0$. $\,\,\Box $

\smallskip
\smallskip
As a consequence of the above theorem, we have the following
interesting special cases which retrieve some results of \cite{r78},
\cite{r3}, and \cite{r62}:\newline

\noindent \textbf{Theorem C.}
\newline
\emph{{\em(i)} Let the $\beta $-conformal change $L\longrightarrow \,\overline{L}%
=e^{\sigma }L+\beta $ be conformal $(\beta =0)$, then $D_{ij}^{r}$
vanishes identically if and only if $\sigma $~is homothetic.\newline
{\em(ii)} Let the $\beta $-conformal change $L\longrightarrow \,\overline{L%
}=e^{\sigma }L+\beta $ be a $\beta$-change $(\sigma =0)$, then $\,
D_{ij}^{r}$ vanishes identically if and only if $\,b_{i}$ is
Cartan-parallel. \vspace{0.6 cm} \newline}

\noindent \textbf {Acknowledgement}: The author wishes to express
his thanks to his colleague prof. Dr. Nabil L Youssef  for fruitful
discussions, critical reading, comments and the revision of the
manuscript.

\bibliographystyle{amsplain}
\bibliography{ref}

\end{document}